\documentclass[letterpaper, 10 pt, conference]{ieeeconf}  
\usepackage{color}
\usepackage{graphicx}
\usepackage{textcomp}
\usepackage[mathcal]{euscript}
\usepackage{mathrsfs}
\usepackage{amsfonts,amssymb,amsbsy}
\usepackage{graphicx,empheq}
\usepackage[export]{adjustbox}

\newcommand{\mc}{\mathcal}

\renewcommand{\theenumii}{\roman{enumii}}
\renewcommand{\theenumiii}{\alph{enumiii}}
\renewcommand{\thepage}{}

\def\e1{{\varepsilon_{11}}}
\def\b1{{\beta_{11}}}
\def\bp3{{\beta_{33}}}
\def\ep3{{\varepsilon_{33}}}

\newtheorem{thm}{Theorem}[section]

\newtheorem{Lem}{Lemma}[section]
\newtheorem{rmk}{Remark}[section]

\IEEEoverridecommandlockouts                              

\title{\LARGE \bf
Output Control of Smart Beams under Uncertain Dynamic Loads through Non-Collocated Sensors and Actuators }

\author{Mark A. Pinsky$^{1}$ and A. \"Ozkan \"Ozer$^{2}$
\thanks{$^{1}$ is a professor and $^{2}$ is a postdoctoal fellow in the Department of Mathematics \& Statistics, University of Nevada-Reno, 1667 N. Virginia St., Reno, Nevada, USA
        {\tt\small aozer@unr.edu \quad pinsky@unr.edu}}%
}

\begin{document}

\maketitle
\thispagestyle{empty}
\pagestyle{empty}

\begin{abstract}
A problem of vibration control of smart beams was addressed in various publications which primarily utilize collocated sensors and actuators and neglect the effect of measurement noise in the observer design. This paper develops a natural design of an output controller which utilizes an eigenfunction approximation of initial continuous model, eliminates control spillover, and consequently leads to an efficient controller which marginalizes effect of bounded system and measurement disturbances while reducing beam vibrations. It is demonstrated that  this control approach can be attained by a non-collocated actuator and a point-sensor of velocity located nearly anywhere on the beam. We show in simulations that the proposed methodology leads to an efficient reduction of beam vibrations enforced by unknown bounded disturbances.

\end{abstract}
\section{INTRODUCTION}
A problem of vibration control of smart piezoelectric beams and other smart structures has been studied extensively in the past decades; see, for example, review papers \cite{Korkmaz}-\cite{Hurlebaus} in the engineering literature as well as references \cite{Banks-Smith}-\cite{Tucsnak} in the distributed parameter systems literature. This problem is closely connected to a more general problem in controlling dynamics of beams and other elastic structures which was addressed in a number of recent publications, for example, in \cite{Fard}-\cite{Vandegriff}. Different approaches to design of observers for linear PDE-systems were surveyed in \cite{Hidayat}. Observer based output controllers of flexible structures were developed in a number of studies, see \cite{Charravarthy}-\cite{Krstic}, and additional references therein.

As is known there two major thrusts for design of estimators and controllers for PDEs. One - approximates PDEs by a finite-dimensional ODE-systems while the other analyses the entire infinite dimensional systems. The former approach leads to utilizing  of standard design techniques readily available for LTI-systems. However, it was shown in \cite{Balas} that this approach is subject to observation and controller spillovers which may destabilize the corresponding infinite dimensional system or degrade performance of an observer.  This and subsequent studies \cite{Chen}-\cite{Hagen} revealed some conditions under which stability of reduced system implies stability of initial PDEs. Another mentioned drawback of this approach that it sometimes renders unnatural feedback laws which are difficult to interpret \cite{Paranjape}.

Disturbance decoupled observer utilizing a linear feedback control has been developed for a special class of second order PDEs in \cite{Demetriou}, which also has additional references on this subject. This design, however, requires knowledge of the operator associated with spatial distribution of the disturbances.

High-gain observers have been used extensively in linear and nonlinear systems for reducing the effect of bounded disturbances on the estimation errors, \cite{Prasov}. However, the efficacy of these observers is limited by accounting of measurement noise. Different approaches have been proposed to mitigate the effect of noise on the observer performance. These approaches were primarily designed for finite-dimensional systems of ODEs, see \cite{Khalil}-\cite{Bo} and additional references therein.

This paper develops a natural approach to output control of PDEs which utilizes finite dimensional approximation of these models. This technique voids controller spillover, and effectively estimates the effect of measurement noise, spillover  and external disturbances on the performance of the close-loop system. The efficiency of this approach in reducing vibrations of forced smart-beam structures is demonstrated in representative simulations.

\section{Beam-Patch Model}
\begin{figure}[htb]
\centering
\includegraphics[width=3.5in]{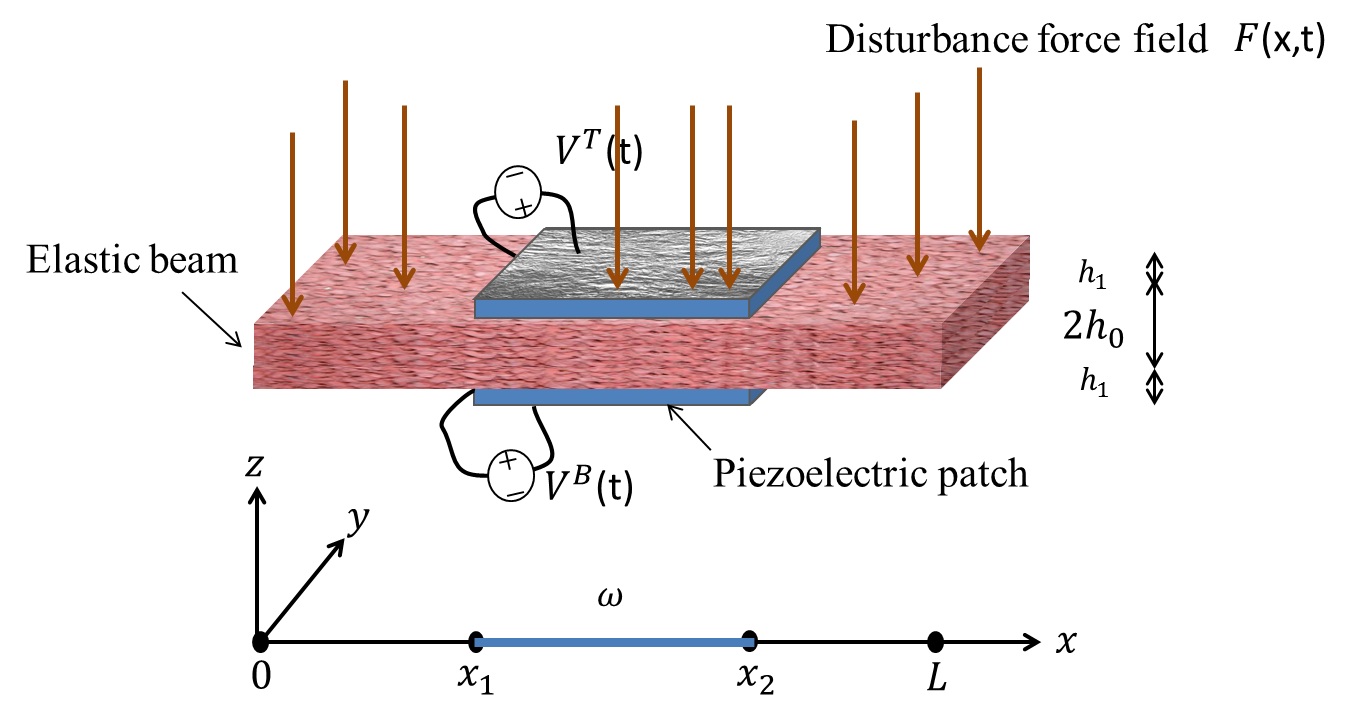}
\caption{\footnotesize An elastic beam with piezoelectric patches which are actuated by voltage $V^T(t)$ at the top and $V^B(t)$ at the bottom.  The disturbance force field $F(x,t)$ is applied in the transverse direction. By choosing $V^B=-V^T$, we force the beam to bend only.}
\label{fig:patch}
\end{figure}
We consider an elastic beam of length $L,$ height $2h,$ and width $b$ with two  piezoelectric patches with height $h_1$ bonded one at the top and another at the bottom of the beam. We also have a sensor located at $0<x_1<L.$ Defining $\omega=[x_1,x_2]$ where  $0<x_1<x_2<L,$ the symmetrically placed patches occupy the region $\Omega_\omega=\omega\times [-b,b]\times \left([-h-h_1, -h]\bigcup [h, h+h_1]\right)\subset \Omega.$  The sensor can also be located on $\omega.$ The patches are insulated at the edges, and no external mechanical stress is applied through the edges. They are also assumed to be bonded perfectly so that no slip occurs. Moreover, each patch is covered with electrodes at lower and upper faces (See Figure \ref{fig:patch}). The model considered in this paper is a modification of the model obtained in \cite{Banks-Smith,accpaper}, where both stretching and bending motions are considered and the only external forces are
 the voltages $V^T(t)$ and $V^B(t)$ which are applied at the top and bottom patches, respectively. In this paper, we consider only the bending motion, i.e.  $V^B=-V^T:=-V,$ with the  disturbance force $F(x,t)$ applied in the transverse direction of the beam. By disregarding the mass and stiffness of the patches, and the magnetic effects (assuming electrostatic or quasi-static electric field), the equation of motion can be written as
\begin{eqnarray}
\nonumber  \rho h \ddot  w  + \frac{E I}{b} w''''-\frac{c_D I}{b}\ddot w' &=&-h\gamma V(t)(\chi_\omega(x))''\\
\label{main} &&+h F(x,t).
 \end{eqnarray}
 where $w, \rho, E, I, c_D,$ and $ \gamma$  denote the transverse displacement, mass density, elasticity constant, mass moment of inertia, damping parameter and  piezoelectric constant. Moreover, $\chi_\omega$ is the characteristic function of the subdomain $w.$ Here  we use the structural damping (square-root damping) in our model to describe the dissipation of the energy of the beam. For this type of damping,  the amplitudes of the normal
modes of vibration are attenuated at rates which are proportional to the oscillation
frequencies.

Dividing (\ref{main}) by $\rho h$  yields
\begin{eqnarray}
\label{main1}    \ddot  w  + \frac{E I}{\rho h b} w''''-\frac{c_D I}{\rho h b}\ddot w' =-\frac{\gamma}{\rho }V(t)(\chi_\omega(x))''+\frac{1}{\rho} F(x,t).
 \end{eqnarray}
Introduction of  dimensionless variables
\begin{eqnarray}
t=\alpha_1 q, \quad x= Lx^*,\quad w=\alpha_4 w^*,
\end{eqnarray}
reduces  (\ref{main1})  to
\begin{eqnarray}
\nonumber    \ddot  w  + \frac{\alpha_1^2 EI}{\rho L^4 h b} w''''-\frac{c_D \alpha_1 I}{\rho h \alpha_1 L^2 b}\ddot w' &=&-\frac{\gamma \alpha_1^2 V(t)(\chi_\omega(x))''}{\rho \alpha_4}\\
\label{main3} &&+\frac{\alpha_1^2}{\rho\alpha_4} F(x,t).
 \end{eqnarray}
 Note that, for simplicity, we use the same notations  for all involved variables.
We define $\alpha_1$ and $\alpha_2$ from the following relations:
\begin{eqnarray}
\nonumber  \frac{\alpha_1^2 EI}{\rho L^4 h b}=1, \quad \frac{\gamma \alpha_1^2 }{\rho \alpha_4}=1
\end{eqnarray}
so that
\begin{eqnarray}
\nonumber  \alpha_1^2= \frac{\rho L^4 h b}{EI}, \quad \alpha_4=\frac{\gamma L^4 hb }{2EI }.
\end{eqnarray}
By letting
\begin{eqnarray}
\nonumber  \frac{c_D \alpha_1 I}{\rho h \alpha_1 L^2 b}=\frac{c_D}{\sqrt{EI\rho h b }}:=a_1,\\
\nonumber \frac{\alpha_1^2}{\alpha_4 \rho}:=a_2,
\end{eqnarray}
write (\ref{main3}) in the dimensionless form
\begin{eqnarray}
\label{nondim}    \ddot  w  +  w''''-a_1\ddot w' =-V(t)(\chi_\omega(x))''+a_2 F(x,t).
 \end{eqnarray}
The hinged boundary  and initial conditions for this equation are
 \begin{eqnarray}
\nonumber w(0,t)=w(1,t)=w''(0,t)=w''(1,t)=0\\
\label{ibc} w(x,0)=w^0,\quad \dot w (x,0)=w^1.
 \end{eqnarray}

\section{Well-posedness}
 First, we consider the homogenous system, i.e., $V\equiv 0$ and $F(x,t)\equiv 0$ in (\ref{nondim}).
Define $z=[w, \dot w]^{\text T}.$ Then (\ref{nondim}) can be written in the state-space form as
\begin{eqnarray}\label{homog}\dot z=\mc A z,\quad z(x,0)=z^0=(w^0,w^1)^{\text T}\end{eqnarray}
where $\mc A=\left(
                 \begin{array}{cc}
                      0 & I\\
                   -D_x^4 & a_1 D_x^2 \\
                 \end{array}
               \right)$
               with ${\text Dom}(\mc A)=H^4_*\times \left(H^2(0,1)\cap H^1_0(0,1)\right)$ where
               $$H^4_*(0,1)=\{w\in H^4(0,1)\cap H^1_0(0,L)~:~ w''(0)=w''(1)=0\}.$$

Let $<u,v>=\int_0^1 u \cdot\bar v~dx$ where $u$ and $v$ be scalar valued functions. Define the bilinear forms
$$c(u,\hat u)=<u'',v''>,\quad a(u,v)=<u,v>.$$
The natural energy associated with (\ref{homog}) is
$$E(t)=\frac{1}{2}\left(a(\dot w,\dot w)+c(w,w)\right), \quad \forall t\in \mathbb{R}^+.$$
Define the Hilbert space $$\mc H=\left(H^2(0,1)\cap H^1_0(0,1)\right)\times H^1_0(0,1)$$
with the energy inner product
$$\left<U,\hat U\right>_{\mc E}:=c(u_1,\hat u_1)+a(u_2,\hat u_2).$$
\begin{Lem} The adjoint of $\mc A$ is
$$\mc A^*=-\mc A(-a_1).$$
Consequently, ${\rm Dom}(\mc A^*)={\rm Dom}(\mc A). $ Moreover, both $\mc A$ and $\mc A^*$ are densely defined and dissipative on $\mc H.$
\end{Lem}
\textbf{Proof:} Let $U,V\in  {\rm Dom} (\mc A).$ Then
\begin{eqnarray}
\nonumber \left<\mc A U,V\right>_{\mc E}&=&c(u_1, v_1)+a(u_2,v_2)\\
\nonumber &=&\int_0^1 \left(u_2'' \bar v_1'' + \left(-u_1'''' +a_1 u_2''\right) \bar v_2\right)~dx\\
\nonumber &=&\int_0^1 \left(u_2'' \bar v_1'' -u_1''\bar v_2''  +a_1 u_2 \bar v_2''\right)~dx\\
\nonumber &=& \left<U, -\mc A(-a_1)V\right>_{\mc E}.
\end{eqnarray}
Then the the conclusion of the first part  of the lemma follows. Obviously, ${\rm Dom}(\mc A)={\rm Dom}(\mc A^*)$ are densely defined in $\mc H.$ Moreover, a simple calculations shows that
$${\rm Re}\left<\mc A U, U\right>_{\mc E}={\rm Re}\left<\mc A^* U, U\right>_{\mc E}=-a_1 \|u_2'\|_{L^2(0,1)}^2.$$
\begin{thm}$\mc A$ is the generator of a $C_0-$semigroup of contractions on $\mc H.$
\end{thm}
\textbf{Proof:} The operator $\mc A$ is closable on $\mc H.$ The conclusion follows from L\"umer Phillips theorem (see \cite{Pazy}). $\square$

\begin{thm} Let $T>0$ and $a_1>0.$ For $V(t)\in L^2(0,T),$ and $F(x,t)\in L^2(0,L;\mathbb{R}^+),$ the solutions of (\ref{nondim})-(\ref{ibc}) satisfy
$w\in (C([0,T];H^2(0,L)\cap H^1_0(0,L))$  and $\dot w\in C([0,T];L^2(0,L)).$
\end{thm}
\textbf{Proof:} Let $y=[w,\dot w]^{\rm T}.$ The system (\ref{nondim})-(\ref{ibc}) can be written as
\begin{eqnarray}\label{mainn}\dot y=\mc A y + B V + a_2 F(x,t)
\end{eqnarray}
where $$B:\mathbb{C} \to \mc H^*$$ defined by $BV=-\left(\begin{array}{c}
0 \\
\vspace{0.1in}
                   \delta'(x_2)-\delta'(x_1)  \\
                  \end{array}\right)V,$
 and $\mc H^*$ is the dual of $\mc H$ with respect to $L^2(0,L).$ The adjoint operator $$B^*: \mc H \to \mathbb{C}$$ defined by $B^*y=   \dot w'(x_2)-\dot w'(x_1).$

 Let $\sigma_n=n\pi, n\in\mathbb{Z}^*:=\mathbb{Z}\backslash \{0\}.$ A short calculation shows that the operator $\mc A$ has the eigenpairs
\begin{eqnarray}\label{eigs}\lambda_n=\sigma_n^2 a=\sigma_n^2 \left\{
  \begin{array}{ll}
    \frac{-a_1\mp \sqrt{a_1^2-4}}{2}, & a_1>2 \\
    \frac{-a_1\pm i \sqrt{4-a_1^2}}{2}, & a_1<2
  \end{array}
\right.
\end{eqnarray}
and
\begin{eqnarray}_n=\frac{\sqrt{2}a}{\sqrt{a^2+1}}\left(\begin{array}{c}
\frac{1}{\lambda_n}\sin{\sigma_n x} \\
\vspace{0.1in}
                   \sin{\sigma_n x}  \\
                  \end{array}\right).
                  \end{eqnarray}

 As $V\equiv 0,$ $F(x,t)\equiv 0,$ and $z^0\in\mc H,$ the equation (\ref{homog}) has a unique $L^2-$analytic solution satisfying $w\in C([0,T];H^2(0,L)\cap H^1_0(0,L))$ and $ \dot w\in C([0,T];L^2(0,L)) $ Moreover, $\mc A$ generates an exponentially stable semigroup on $\mc H,$ see i.e. \cite{Trig},\cite{Las}.

The solution of (\ref{homog}) is $y(x,t)=\sum_{n\in \mathbb Z^*} c_n e^{i\lambda_n t} \psi_n$ with $ \{c_n\}_{n\in\mathbb{Z}^*} \in  l_2,$
and $\left\|y^0 \right\|_{\mc E}^2=\sum_{n\in \mathbb Z^*}  |c_n|^2.$

Therefore, there exists a constant $C>0$ such that
\begin{eqnarray}
\nonumber&& \int_0^T \left|B^*e^{A^*t}z\right|^2 ~dt = \int_0^T\left|\dot z'(x_2)-\dot z'(x_1)\right|^2~dt\\
\nonumber && =\int_0^T\left|\sum_{n\in \mathbb Z^*}  i c_n  \lambda_n e^{i\lambda_n t} \left(\phi_n'(x_2) -\phi_n'(x_1)\right)\right|^2~dt\\
 \nonumber &&\le C \sum_{n\in \mathbb Z^*}  |c_n|^2\\
\nonumber  && =C \left\|z^0\right\|_{\mc E}^2
\end{eqnarray}
This implies that $B^*$ is an admissible control operator for (\ref{homog}), and therefore $B$ is an admissible control operator for (\ref{mainn}) with $F(x,t)\equiv 0.$ Since $V(t)\in L^2(0,T),$ and $F(x,t)\in L^2(0,L;\mathbb{R}^+),$ the statement of the theorem follows.
\section{Discrete equations}
 The normalized eigenfunctions for (\ref{nondim})-(\ref{ibc}) are $\psi_n(x)=\{\sqrt{2}\sin{\sigma_n x}\}_{n=1}^{\infty}.$ Now we choose $$w(x,t)=\sum_{n=1}^{N}w_n(t)\psi_n(x), \quad F(x,t)=\sum_{n=1}^{N}f_n(t)\psi_n(x),$$
 multiply the equation (\ref{nondim}) by $\psi_n,$ and integrate by parts to get
\begin{eqnarray}
\label{discrete}    \ddot  w_n + \sigma_n^4 w_n  +  a_1\sigma_n^2 \dot w_n  =V(\psi_n'(x_2)-\psi_n'(x_1))+a_2 f_n
 \end{eqnarray}
for $n=1,2,\cdots, N.$ Let $$z_N=[w_1, w_2, \ldots w_N, \dot w_1, \dot w_2, \ldots, \dot w_N]^T.$$  We formulate (\ref{discrete}) as the following:
\begin{eqnarray}
\label{dum}\dot z_N= \mc A_N z_N + B_N V(t)+F_N, \quad z(x,0)=z^0
\end{eqnarray} where
 $\mc A_N=\left(
                 \begin{array}{cc}
                      0 & I_N\\
                   -A_N & -Q_N \\
                 \end{array}
               \right),$ and
                             \begin{align}\nonumber & I_N={\rm diag}(1, \cdots, 1), \\
\nonumber              & A_N={\rm diag}(\sigma_1^4, \cdots, \sigma_N^4),\\
           \nonumber   & Q_N={\rm diag}(a_1\sigma_1^2, \cdots, a_1\sigma_N^2),\\
           \nonumber & B_N=[0,\cdots,0,\psi_1(x_2) -\psi_1(x_1),\cdots, \psi_N(x_2) -\psi_N(x_1)]^{\rm T},\\
           & F_N=[ 0,0,\cdots, 0, f_1,\cdots f_N]^{\rm T}.
               \end{align}

\subsection{Controller Design}
We choose an integral-type feedback controller in the form
\begin{eqnarray}
\label{feedback}V(t)=\int_0^1 k(x) \dot w(x,t)~ dx
\end{eqnarray}
where the gain function $k(x)\in L^2(0,1).$  Since $\{\psi_n\}_{n=1}^{\infty}$ is a basis for $L^2(0,1),$ we can write
\begin{eqnarray}V(t)&=&\int_0^L \left(\sum_{n=1}^{\infty} k_n\psi_n(x)\right)\left(\sum_{n=1}^{\infty} \dot w_n \psi_n(x)\right)~dx \quad~~\\
\label{oz1} &=&\sum_{n=1}^{\infty} k_n\dot w_n.
\end{eqnarray}
Now we choose $k_n=0$ for all $n>N$ such that (\ref{oz1}) can be written as
\begin{eqnarray}V(t)&=&\int_0^L \left(\sum_{n=1}^{N} k_n\psi_n(x)\right)\left(\sum_{n=1}^{\infty} \dot w_n \psi_n(x)\right)~dx \quad~~\\
\label{oz2} &=&\sum_{n=1}^{N} k_n\dot w_n.
\end{eqnarray}
\subsection{Observers}
Let $x_0\in (0,1)$ and the output be
\begin{eqnarray} y(t)&=&s_1  w(x_0,t)+s_2 \dot w(x_0,t)+\xi(t)\\
\nonumber &=& s_1 \sum_{n=1}^{N} w_n(t)\psi_n(x_0)\\
&&\quad\quad+ s_2 \sum_{n=1}^{N} \dot w_n(t)\psi_n(x_0)+ r(t) + \xi(t).~~
\end{eqnarray}
where $\xi(t)$ is the bounded measurement noise, i.e. $|\xi|<\infty,$  and $r(t)$ is the residual of the truncated eigenfunction expansion, and $s_1, s_2\in \mathbb{R}$ are adjustable weights.
Then the observation matrix is
$$C_N=[s_1 \psi_1, \cdots s_1\psi_N(x_0), s_2 \psi_1(x_0),\cdots, \psi_N(x_0)]^{\rm T}.$$
In this paper we assume that the velocity of the beam can be measured at a certain point on the beam. So $s_1=0$ and $s_2=1.$ Now the combined system is
\begin{eqnarray}
\label{dum2}\dot z_N= \mc A_N z_N + B_N V(t)+F_N(t), \quad z(x,0)=z^0\\
\label{dum3}y_N=C_N z_N + \varepsilon ; \quad \varepsilon=r(t)+\xi (t)
\end{eqnarray}  Note that $Cz=\dot w(x_0,t)$ is an admissible observation operator for the homogenous system (\ref{homog}); i.e. $ \int_0^T \|Cz\|^2 ~dt \le \sum_{n\in Z^*} c_n^2 \le \|z^0\|_{\varepsilon}^2.$ Therefore the boundedness of $r(t)=\sum_{n=N+1}^{\infty} C_i z_i,$ and therefore $y(t),$ follows by the classical perturbation argument.
\section{Observer Design}
Choosing $V=-K_N \hat z_N,$ we write the observer equation as
\begin{eqnarray}
\nonumber \dot {\hat z}_N=\left(\mc A_N  - B_N K_N\right) {\hat z}_N + L_N(y_N-C_N \hat z_N)
\end{eqnarray}
where $L_N$ and $K_N$ are the gain matrices for observer and controller, respectively. Next, setting the estimation error $e_N:=z_N-\hat z_N$ yields the coupled state and error dynamics
\begin{eqnarray}
\nonumber \dot { z}_N= \left(\mc A_N  - B_N K_N\right){z}_N + B_N K_N e_N + F_N\\
\label{ozzz} \dot e_N=(\mc A_N  -L_NC_N )e_N-L_N\varepsilon + F_N.
\end{eqnarray}
Due to the separation principle, eigenvalues of the matrices $\mc A_N  - B_N K_N$  and $\mc A_N  - L_N C_N$ can be set up to any
given values by an appropriate choice of $K_N$ and $L_N$ if the pairs $(A_N, B_N)$ and $(A_N,C_N)$ are controllable and observable, respectively. In our case, we choose $K_N$ and $L_N$ to marginalize the effect of disturbances and the estimation error on the system dynamics.

\begin{thm}
The pair $(A_N,C_N)$ is observable if $x_0\ne \frac{k}{n}$ where $k\le n,$ $n\le N,$ $n\in \mathbb{Z}^+,$ and $k\in \mathbb{Z}^+.$ Moreover, the pair $(A_N,B_N)$ is controllable if $x_2\ne x_1+2 n k,$ where $0\le x_1<x_2\le 1,$ $N\ge k\ge 2n,$ $n\in \mathbb{Z}^+,$ and  $k\in \mathbb{Z}^+.$
\end{thm}

\textbf{Proof:} Obviously, the matrix $A_N$ is diagonalizable. Since every single entry of the matrices $B_N$  and $C_N$ are nonzero due to the respective conditions for each matrix, the statement of the theorem follows.

 \subsection{Error Analysis}
For the rest of the paper, $\|\cdot \|$ refers to the Euclidean norm. Accounting for the external forces, the measurement noise and spillover in observations jeopardize the convergency of $\|e\|$ and $\|z\|$ to zero. However, these norms can be minimized by the appropriate choice of $L_N$ and $K_N.$
Let $A_L=A_N-L_NC_N$ and $A_K=A_N-B_NK_N.$ The norm of the solutions to (\ref{ozzz}) can be bounded as follows
\begin{eqnarray}
\label{dumm1}\|e_N\|&\le& e^{-\lambda_L t} e(0) + \frac{\|F(t)\| + \|L_N\|\|\varepsilon\|}{\lambda_L},\\
\label{dumm2}\|z_N\|&\le& e^{-\lambda_K t} z(0) + \frac{\|F(t)\| + \|B_NK_N\|\|e\|}{\lambda_K}\quad
\end{eqnarray}
where $\lambda_L , \lambda_K>0$ are the minimums of the absolute values of the real part of the eigenvalues of the matrices $A_L$ and $A_K,$ respectively. Both bounds include steady-state and transient components, and the later decay exponentially to zero. Observe that the bound for the steady-state error term in (\ref{dumm1}), which is  called $\|e_{\rm st}\|,$  does not approach zero as $\lambda_L\to \infty$ since in this case $\|L_N\|\to \infty.$ Consequently, a practical solution is to find  an $L_N$ which minimizes $\|e_{\rm st}\|.$ This implies that $\|e_N\|<\infty,$ then $K_N$ is chosen  such  that the norm of the steady-steady component $\|z_{st}\|$ in (\ref{dumm2}) is minimized likewise. Note that, in the practical applications, it is commonly assumed that $\lambda_K<\lambda_L.$

\subsection{Residual bounds}
The contribution of the uncontrolled modes, i.e. $\sum_{i=N+1}^{\infty} w_i(t) \psi_i(x),$ is known as the residual or spillover effect.
Let us now estimate the norm of solutions of uncontrolled infinite-dimensional residual system
\begin{eqnarray}
\label{res1} \ddot w_k + a_1 \sigma_k^2 \dot w_k + \sigma_k^4 w_k=a_2 f_k, \quad k\ge N+1.
\end{eqnarray}
These mutually decoupled system of  equations  are also decoupled from the system (\ref{dum2})-(\ref{dum3}) due to the choice of the controller (\ref{oz2}).  Obviously, by (\ref{eigs}), the real parts of eigenvalues of this system are negative and bounded from above, i.e.
\begin{eqnarray}
\label{res2} {\rm Real} (\lambda_k) \le -a_1 \sigma_{k}^2=-a_1 \pi^2 k^2; \quad k\ge N+1    .
\end{eqnarray}
Thus
\begin{eqnarray}
\label{res3} \|w_k\|\le \frac{a_2 \|f_k\|}{a_1\pi^2 k^2}; \quad k\ge N+1.
\end{eqnarray}
Let us assume that
\begin{eqnarray}
\label{res4} \|f_k\|\le f^0, \quad k\ge N+1.
\end{eqnarray}
Then
\begin{eqnarray}
\nonumber \sum_{k=N+1}^{\infty} \|w_k\| &\le& \frac{a_2 f^0}{a_1 \pi^2} \left(\sum_{k=N+1}^{\infty} \frac{1}{k^2}\right)\\
\label{res5} &\le & \frac{a_2 f^0}{a_1 \pi^2} \frac{1}{N+1}
\end{eqnarray}
Thus under (\ref{res4}), the norm of  spillover component in the solution decays as $O\left(\frac{1}{N+1}\right)$. In the same time, a common assumption that $F(x,t)$  is continuously differentiable in $x$ yields that $f_k=O\left(\frac{1}{k^2}\right)$ which enhances residual estimate as follows
\begin{eqnarray}
\label{res5}  \sum_{k=N+1}^{\infty} \|w_k\| &\le & \frac{a_2 f^0}{3a_1 \pi^2} \frac{1}{(N+1)^3}.
\end{eqnarray}
This rapid convergence to zero justifies the use of relatively low-dimensional models in control applications of continuous systems.

\begin{rmk} In this paper we assume that the structural damping in our model, i.e. the term $a_1\dot w''.$ Note that the alteration of this damping with the Kelvin-Voigt damping, i.e. $a_1\dot w'''',$ amplifies the contribution of residual modes since in this case (\ref{res2}) is replaced by
\begin{eqnarray}
\label{res6} {\rm Real} (\lambda_k) \le C.
\end{eqnarray}
where $C<0$ is independent from $k.$
\end{rmk}

\section{Simulations}
\begin{figure}
\centering
\includegraphics[width=8.8cm,height=3.5cm,center]{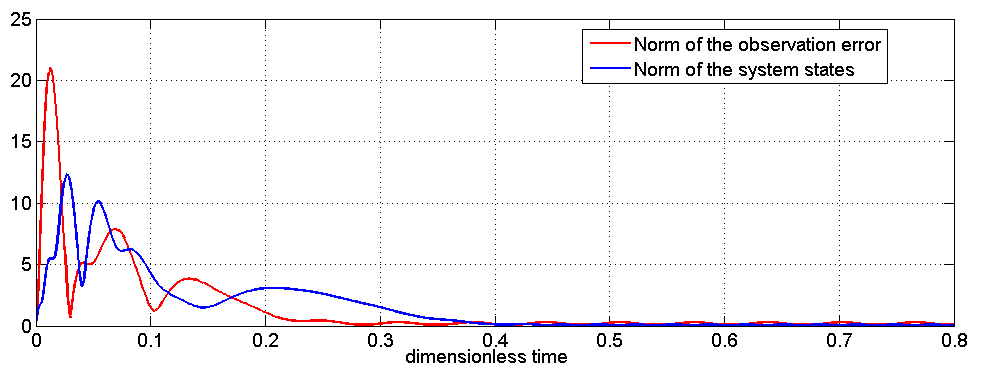}
\vspace{-0.2in}
\caption{Time-histories of the norms of observation errors and system-states of the finite-dimensional system with $\lambda_L=34.$}
\label{fig1}
\end{figure}

\begin{figure}[]
\centering
\includegraphics[width=8.8cm,height=3.5cm,center]{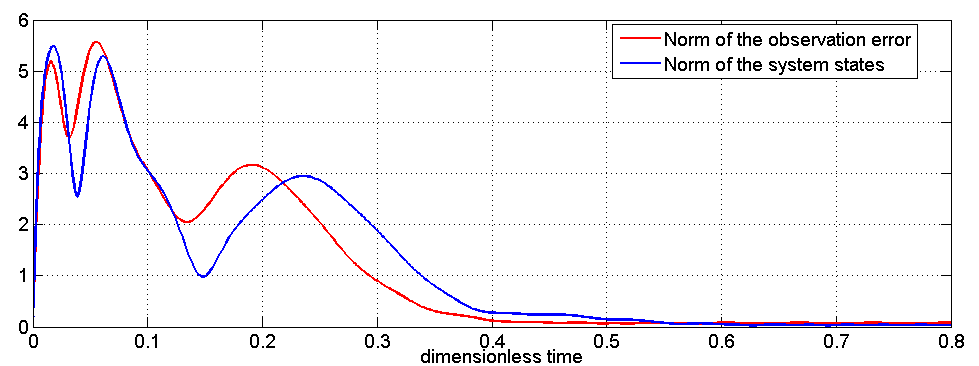}
\vspace{-0.2in}
\caption{Picking phenomenon is more pronounced for larger observer gain with $\lambda_L=64.$ }
\label{fig2}
\end{figure}

\begin{figure}
\centering
\includegraphics[width=8.8cm,height=3.5cm,center]{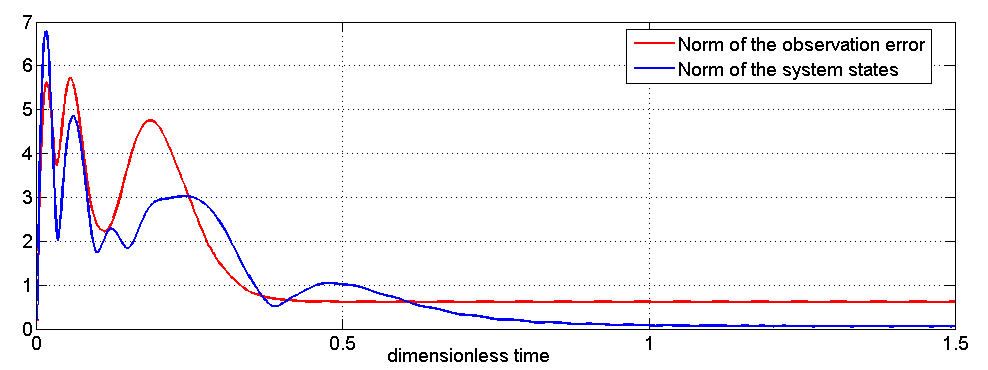}
\vspace{-0.2in}
\caption{The symmetrical patches are located at $\omega=(0,0.1)$, the velocity sensor is located at the left edge of the patched region $x_0=0.095,$ and $N=3.$}
\label{fig3}
\end{figure}

\begin{figure}
\centering
\includegraphics[width=8.8cm,height=3.5cm,center]{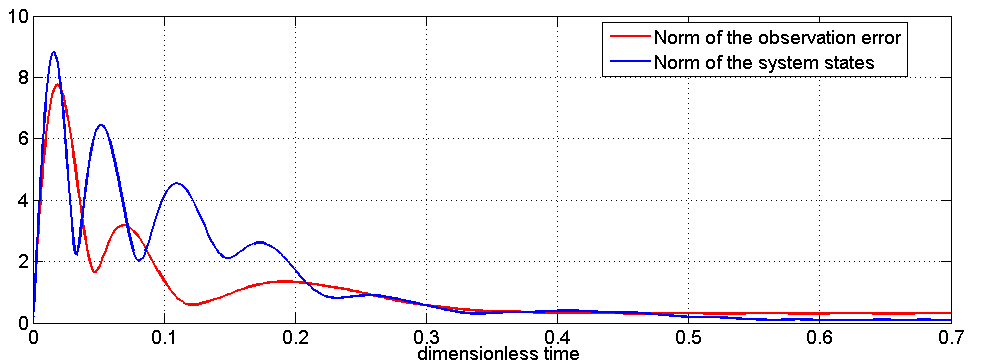}
\caption{The symmetrical patches are located at $\omega=(0,0.1)$, the velocity sensor is located at the right edge of the patched region $x_0=0.98,$ and $N=3.$}
\label{fig4}
\end{figure}

\begin{figure}
\centering
\includegraphics[width=8.8cm,height=3.5cm,center]{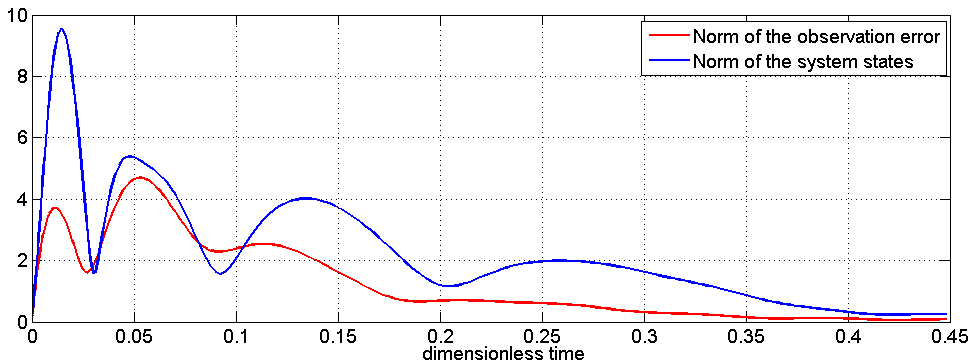}
    \caption{The symmetrical patches are located at $\omega=(0,0.1)$,  the velocity sensor is located close to the middle of the beam $x_0=0.6,$ and $N=3.$}
    \label{fig5}
\end{figure}

\begin{figure}
\centering
\includegraphics[width=9cm,height=3.5cm,center]{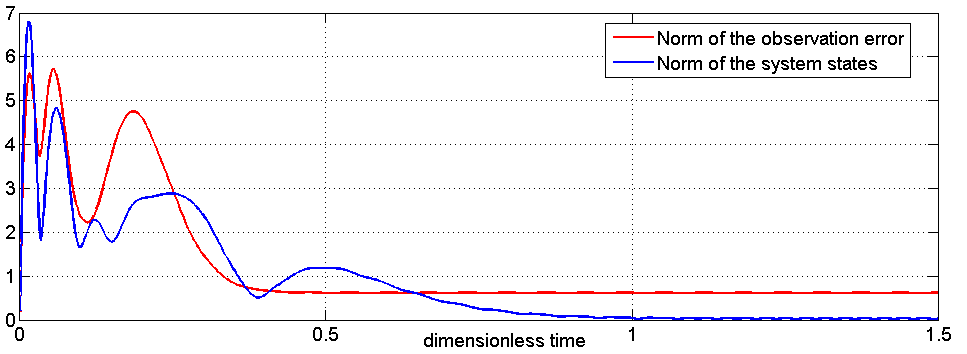}
    \caption{The lengths of patches are  smaller and are located  at $\omega=(0,10^{-8})$, the velocity sensor is located close to the middle of the beam $x_0=0.6,$ and $N=3.$}
    \label{fig6}
\end{figure}

\begin{figure}
\centering
\includegraphics[width=9.2cm,height=3.5cm, center]{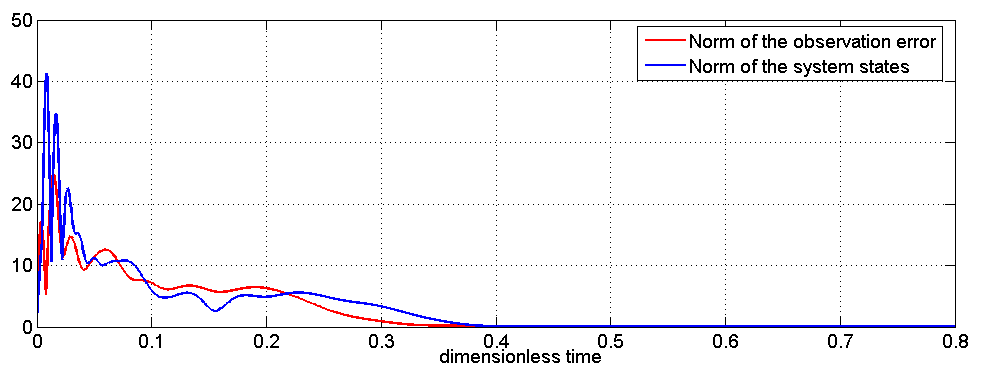}
\vspace{-0.2in}
    \caption{ The symmetrical patches are located at $\omega=(0,0.1)$, the velocity sensor is located close to the right edge of the beam $x_0=0.98,$ and $N=5.$}
    \label{fig7}
    \vspace{-0.2in}
\end{figure}

To demonstrate the efficiency of the proposed design of output controller, we present time-histories of the norms of  error and state vector in Fig. \ref{fig1}-\ref{fig7}. We choose $N=3$  for Fig. \ref{fig1}-\ref{fig6}, and $N=5$ for Fig. \ref{fig7}, $c_D=0.01$. The force field is modeled by equal polyharmonic forces acting on the first three vibration modes. Each forcing term has 11 harmonic components including a resonance component with the  lowest vibration frequency. The maximum value of each forcing term is bounded by 11 in our simulations.

Comparison of Fig. \ref{fig1} and Fig. \ref{fig2} demonstrates that the so-known picking phenomenon becomes more pronounced for amplified observer gains while these gains condense the durations of the transient behavior. In these two cases the position of the sensor is collocated at the edge of the actuator patch.

 Figs. \ref{fig3}-\ref{fig5} show time-histories of error and state-vector norms in the cases when the sensor is located near the left or right edges, or the middle point of the bar while the actuator patch is extended from the left edge of the bar, i.e., $x_1=0, x_2=0.1.$  Fig. \ref{fig6} shows the reduction of the steady-state solution obtained for very short patches, i.e. the $x_2-x_1=10^{-8}.$ Moreover, Fig \ref{fig7} shows the dynamics of the system  (\ref{ozzz}) in case of $N=5.$


 Note that in all cases mentioned above the superior reduction of steady-state vibrations is developed by the proposed control. In fact, the norm of external forces exceeds the norm of steady- state solutions in all cases about two order in magnitude despite the fact that these forces resonate with lowest frequency and the damping coefficient is assumed to be small.

\section{Discussions and Conclusions}
This paper addresses the problem of designing an efficient output controller that leads to a significant reduction of vibrations in an Euler-Bernoulli beam subject to unknown but bounded force field. While this approach is based on the finite-dimensional approximation of the corresponding partial differential equation model (\ref{main}), it naturally avoids the spillover caused by control, and leads to an efficient treatment of the measurement spillover.			
We establish observability and controllability conditions which specify exceptional locations of a sensor/an actuator violating one/both of these conditions. Thus it is demonstrated that the non-collocated sensor and actuators can be placed almost anywhere on the bar to attain quite similar vibration reduction.

Our design accounts for the bounded observation noise which elevates the norms of observation error. We bound the norms of observation errors and system states. We argue that the observer and controller gains can be chosen to substantially reduce these norms. 							

Next, we bound the norm of uncontrolled residual components in this finite-dimensional approximation and estimate the speed with which this norm approaches zero as the number of controlled modes increases. This leads to a sound assessment of the accuracy and required dimensions of such approximations. Finally, it follows from the presented simulations that the proposed control architecture facilitates significant reduction of beams vibrations for almost arbitrary positions of a sensor/ an actuator.

In this paper a reduction of the steady state vibrations of a beam is attained due to application of high-gain observers and controllers which, in turn, develops a transient behavior known as the picking phenomenon \cite{Prasov} which can be mitigated in practical applications by saturating  the controller \cite{Esfan} or by delaying control engagement until the observer error falls under a certain threshold, or by application of more adaptive piece-wise linear observer designs \cite{Prasov}.  While relatively straightforward in this case, the implementation of these approaches is left aside in this paper.

While in this paper the developed methodology is applied to a relatively simple model of Euler-Bernoulli beam, it has potential to be useful in similar applications to more accurate models as  Mindlin –-Timoshenko beam (\cite{accpaper},\cite{Dietl})  and to more complex elastic structures as plates, shells, etc.  A more accurate modeling of electric and magnetic properties of piezoelectric actuators \cite{accpaper} and accounting for nonlinear components should lead to more comprehensive and precise modeling of smart structures which will be the topic of our future research.

\end{document}